# ANTICIPATED BACKWARD STOCHASTIC DIFFERENTIAL EQUATIONS[1]

By Shige Peng and Zhe Yang

*Shandong University, and Shandong University and Cambridge University*

In this paper we discuss new types of differential equations which we call anticipated backward stochastic differential equations (anticipated BSDEs). In these equations the generator includes not only the values of solutions of the present but also the future. We show that these anticipated BSDEs have unique solutions, a comparison theorem for their solutions, and a duality between them and stochastic differential delay equations.

**1. Introduction.** Consider these types of stochastic differential delay equations (SDDEs):

$$(1)\qquad \begin{cases} dX_t = (\mu_t X_t + \bar{\mu}_{t-\theta} X_{t-\theta})\,dt + (X_t \sigma_t^T + X_{t-\theta}\bar{\sigma}_{t-\theta}^T)\,dW_t, \\ \qquad\qquad t \in [t_0, T+\theta]; \\ X_t = x_t, \qquad t \in [t_0 - \theta, t_0], \end{cases}$$

where $W$ is a $d$-dimensional Brownian motion, $\theta > 0, x_t$ is a deterministic function, and $Q$ is a given $\mathscr{F}_T^W$-measurable random variable. In the case where $\bar{\mu} = \bar{\sigma} \equiv 0$, this model is very typical in finance as the price of a stock. Then $Y_{t_0} = E[X_T Q | \mathscr{F}_{t_0}]$ can be the price of an option valued $Q$ at maturity time $T$ if $x_t \equiv 1$. It is easy to prove that (see, e.g., El Karoui, Peng and Quenez [7]) $Y_\cdot$ is a solution to the following backward stochastic differential equation (BSDE):

$$-dY_t = (\mu_t Y_t + Z_t \sigma_t)\,dt - Z_t\,dW_t, \qquad Y_T = Q.$$

This SDE with delay, in which $\bar{\mu}$ and $\bar{\sigma}$ are nonzero, has a solution. An interesting question is whether it can be expressed in the form of equation (1).

Received April 2007; revised March 2008.
[1]Supported in part by the National Basic Research Program of China (973 Program) Grant 2007CB814900 (Financial Risk).

The answer is positive if we can solve the following new type of "anticipated" BSDE:

$$
(2)\qquad \begin{cases} -dY_t = (\mu_t Y_t + \bar{\mu}_t E^{\mathscr{F}_t}[Y_{t+\theta}] + Z_t\sigma_t \\ \qquad\qquad + E^{\mathscr{F}_t}[Z_{t+\theta}]\bar{\sigma}_t + l_t)\,dt - Z_t\,dW_t, & t\in[t_0,T]; \\ Y_t = Q_t, & t\in[T,T+\theta]; \\ Z_t = P_t, & t\in[T,T+\theta]. \end{cases}
$$

We observe that the generator, that is, the $dt$ part of the BSDE, contains the values $(Y_\cdot, Z_\cdot)$ for present time $t$ as well as for future time $t+\theta$. This is a new duality phenomenon for SDEs and BSDEs.

In this paper we consider a more general form of this new type of BSDE:

$$
\begin{cases} -dY_t = f(t, Y_t, Z_t, Y_{t+\delta(t)}, Z_{t+\zeta(t)})\,dt - Z_t\,dW_t, & t\in[0,T]; \\ Y_t = \xi_t, & t\in[T,T+K]; \\ Z_t = \eta_t, & t\in[T,T+K]. \end{cases}
$$

The paper is organized as follows. In Section 2 we consider the duality between SDDEs and anticipated BSDEs. After a brief presentation of some known results that we will use in Section 3, we prove an existence and uniqueness result for anticipated BSDEs in Section 4. In Section 5 we give an important result for anticipated BSDEs: a comparison theorem. In Section 6 we use the duality between SDDEs and anticipated BSDEs mentioned in Section 2 to solve a stochastic control problem.

**2. Duality between SDDEs and anticipated BSDEs.** It is well known that there is perfect duality between SDEs and BSDEs (see El Karoui, Peng and Quenez [7]). In this section we consider duality between the SDDEs and the anticipated BSDEs mentioned above. We will use this duality to solve a stochastic control problem in Section 6.

THEOREM 2.1. *Suppose $\theta>0$ is a given constant and $\mu_\cdot, \bar{\mu}_\cdot \in L^2_{\mathscr{F}}(t_0-\theta, T+\theta), l_\cdot \in L^2_{\mathscr{F}}(t_0,T), \sigma_\cdot, \bar{\sigma}_\cdot \in L^2_{\mathscr{F}}(t_0-\theta, T+\theta;\mathbb{R}^{d\times 1}), \mu_\cdot, \bar{\mu}_\cdot, \sigma_\cdot, \bar{\sigma}_\cdot$ are uniformly bounded. Then for all $Q_\cdot \in S^2_{\mathscr{F}}(T,T+\theta), P_\cdot \in L^2_{\mathscr{F}}(T,T+\theta;\mathbb{R}^d)$, the solution $Y_\cdot$ of the anticipated BSDE (2) can be given by the closed formula*

$$
Y_t = E^{\mathscr{F}_t}\left[X_T Q_T + \int_t^T X_s l_s\,ds + \int_T^{T+\theta} (Q_s\bar{\mu}_{s-\theta} + P_s\bar{\sigma}_{s-\theta})X_{s-\theta}\,ds\right], \qquad a.e.,\quad a.s.,
$$

*where $X_s$ is the solution to the SDDE*

$$
(3)\qquad \begin{cases} dX_s = (\mu_s X_s + \bar{\mu}_{s-\theta}X_{s-\theta})\,ds + (X_s\sigma_s^T + X_{s-\theta}\bar{\sigma}_{s-\theta}^T)\,dW_s, \\ \qquad\qquad s\in[t,T+\theta]; \\ X_t = 1, \\ X_s = 0, \qquad s\in[t-\theta,t]. \end{cases}
$$

Proof. First, we show that (3) has a unique solution. When $s \in [t, t+\theta]$, (3) becomes

$$\begin{cases} dX_s = \mu_s X_s\, ds + X_s \sigma_s^T\, dW_s, & s \in [t, t+\theta]; \\ X_t = 1. \end{cases} \tag{4}$$

We can then easily obtain a unique continuous solution $\varsigma_\cdot$ for (4). When $s \in [t+\theta, T+\theta]$, (3) becomes

$$\begin{cases} dX_s = (\mu_s X_s + \bar{\mu}_{s-\theta} X_{s-\theta})\, ds + (X_s \sigma_s^T + X_{s-\theta} \bar{\sigma}_{s-\theta}^T)\, dW_s, \\ \qquad\qquad s \in [t+\theta, T+\theta], \\ X_s = \varsigma_s, \qquad s \in [t, t+\theta]. \end{cases} \tag{5}$$

Equation (5) is a classical SDDE, thus, it has a unique solution. Applying Itô's formula to $X_s Y_s$ for $s \in [t, T]$ and taking conditional expectations under $\mathscr{F}_t$, we get

$$\begin{aligned} &E^{\mathscr{F}_t}[X_T Y_T] - X_t Y_t \\ &\quad = E^{\mathscr{F}_t}\biggl[\int_t^T (Y_s \bar{\mu}_{s-\theta} X_{s-\theta} - E^{\mathscr{F}_s}[Y_{s+\theta}] \bar{\mu}_s X_s \\ &\qquad\qquad + Z_s \bar{\sigma}_{s-\theta} X_{s-\theta} - E^{\mathscr{F}_s}[Z_{s+\theta}] \bar{\sigma}_s X_s - X_s l_s)\, ds\biggr]. \end{aligned}$$

Because $X_t = 1$ and $X_s = 0$, $s \in [t-\theta, t)$, we have

$$\begin{aligned} Y_t &= E^{\mathscr{F}_t}\biggl[X_T Y_T + \int_t^T X_s l_s\, ds\biggr] - E^{\mathscr{F}_t}\biggl[\int_t^T (Y_s \bar{\mu}_{s-\theta} X_{s-\theta} - Y_{s+\theta} \bar{\mu}_s X_s)\, ds\biggr] \\ &\quad - E^{\mathscr{F}_t}\biggl[\int_t^T (Z_s \bar{\sigma}_{s-\theta} X_{s-\theta} - Z_{s+\theta} \bar{\sigma}_s X_s)\, ds\biggr] \\ &= E^{\mathscr{F}_t}\biggl[X_T Y_T + \int_t^T X_s l_s\, ds - \int_t^T Y_s \bar{\mu}_{s-\theta} X_{s-\theta}\, ds + \int_{t+\theta}^{T+\theta} Y_s \bar{\mu}_{s-\theta} X_{s-\theta}\, ds\biggr] \\ &\quad - E^{\mathscr{F}_t}\biggl[\int_t^T Z_s \bar{\sigma}_{s-\theta} X_{s-\theta}\, ds - \int_{t+\theta}^{T+\theta} Z_s \bar{\sigma}_{s-\theta} X_{s-\theta}\, ds\biggr] \\ &= E^{\mathscr{F}_t}\biggl[X_T Q_T + \int_t^T X_s l_s\, ds + \int_T^{T+\theta} (Q_s \bar{\mu}_{s-\theta} X_{s-\theta} + P_s \bar{\sigma}_{s-\theta} X_{s-\theta})\, ds\biggr]. \qquad \square \end{aligned}$$

**3. Preliminaries.** Let $(\Omega, \mathscr{F}, P, \mathscr{F}_t, t \geq 0)$ be a complete stochastic basis such that $\mathscr{F}_0$ contains all $P$-null elements of $\mathscr{F}$ and suppose that the filtration is generated by a $d$-dimensional standard Brownian motion $W = (W_t)_{t \geq 0}$. Given $T > 0$, denote the norm in $\mathbb{R}^m$ by $|\cdot|$. We will use the following notation:

- $L^2(\mathscr{F}_T; \mathbb{R}^m) = \{\mathbb{R}^m$-valued $\mathscr{F}_T$-measurable random variables such that $E[|\xi|^2] < \infty\}$;

- $L^2_{\mathscr{F}}(0,T;\mathbb{R}^m) = \{\mathbb{R}^m$-valued and $\mathscr{F}_t$-adapted stochastic processes such that $E[\int_0^T |\varphi_t|^2\, dt] < \infty\}$;
- $S^2_{\mathscr{F}}(0,T;\mathbb{R}^m) = \{$continuous processes in $L^2_{\mathscr{F}}(0,T;\mathbb{R}^m)$ such that $E[\sup_{0\le t\le T} |\varphi_t|^2] < \infty\}$.

If $m=1$, we denote them by $L^2(\mathscr{F}_T), L^2_{\mathscr{F}}(0,T)$ and $S^2_{\mathscr{F}}(0,T)$. The above $L^2$ are all separable Hilbert spaces.

The following lemmas can be found in Peng [13], Section 3. For their originalities we refer to the notes of [13] or [7]. Our Lemma 3.1 is Lemma 3.1 of Peng [13]. Lemma 3.2, which is Theorem 3.2 of Peng [13], is a basic result of BSDEs: an existence and uniqueness theorem. Both Lemmas 3.3 and 3.4 are comparison theorems for solutions of BSDEs. Lemma 3.3 is Theorem 3.3 of Peng [13] and can also be found in El Karoui, Peng and Quenez [7]. Lemma 3.4 can be easily obtained from Lemma 3.3.

Lemma 3.1. *For a fixed $\xi \in L^2(\mathscr{F}_T)$ and $g_0(\cdot)$ which is an $\mathscr{F}_t$-adapted process satisfying $E[(\int_0^T |g_0(t)|\, dt)^2] < +\infty$, there exists a unique pair of processes $(y_\cdot, z_\cdot) \in L^2_{\mathscr{F}}(0,T;\mathbb{R}^{1+d})$ satisfying the following BSDE:*

$$y_t = \xi + \int_t^T g_0(s)\, ds - \int_t^T z_s\, dW_s, \qquad t\in[0,T].$$

*If $g_0(\cdot) \in L^2_{\mathscr{F}}(0,T)$, then $(y_\cdot, z_\cdot) \in S^2_{\mathscr{F}}(0,T) \times L^2_{\mathscr{F}}(0,T;\mathbb{R}^d)$. We have the following basic estimate:*

$$\begin{aligned} &|y_t|^2 + E^{\mathscr{F}_t}\left[\int_t^T \left(\frac{\beta}{2}|y_s|^2 + |z_s|^2\right) e^{\beta(s-t)}\, ds\right] \\ &\qquad \le E^{\mathscr{F}_t}[|\xi|^2 e^{\beta(T-t)}] + \frac{2}{\beta} E^{\mathscr{F}_t}\left[\int_t^T |g_0(s)|^2 e^{\beta(s-t)}\, ds\right]. \end{aligned} \tag{6}$$

*In particular,*

$$\begin{aligned} &|y_0|^2 + E\left[\int_0^T \left(\frac{\beta}{2}|y_s|^2 + |z_s|^2\right) e^{\beta s}\, ds\right] \\ &\qquad \le E[|\xi|^2 e^{\beta T}] + \frac{2}{\beta} E\left[\int_0^T |g_0(s)|^2 e^{\beta s}\, ds\right], \end{aligned} \tag{7}$$

*where $\beta > 0$ is an arbitrary constant. We also have*

$$E\left[\sup_{0\le t\le T} |y_t|^2\right] \le kE\left[|\xi|^2 + \int_0^T |g_0(s)|^2\, ds\right], \tag{8}$$

*where the constant $k$ depends only on $T$.*

We assume that $g = g(\omega, t, y, z) : \Omega \times [0,T] \times \mathbb{R}^m \times \mathbb{R}^{m\times d} \longrightarrow \mathbb{R}^m$ satisfies the following conditions:

(a) $g(\cdot,y,z)$ is an $\mathbb{R}^m$-valued and $\mathscr{F}_t$-adapted process satisfying the Lipschitz condition in $(y,z)$, that is, there exists $\rho>0$ such that, for each $y,y'\in\mathbb{R}^m$ and $z,z'\in\mathbb{R}^{m\times d}$, $|g(t,y,z)-g(t,y',z')|\le\rho(|y-y'|+|z-z'|)$.

(b) $g(\cdot,0,0)\in L^2_{\mathscr{F}}(0,T;\mathbb{R}^m)$.

Lemma 3.2. *Assume that $g$ satisfies* (a) *and* (b), *then for any given terminal condition $\xi\in L^2(\mathscr{F}_T;\mathbb{R}^m)$, BSDE*

$$Y_t=\xi+\int_t^T g(s,Y_s,Z_s)\,ds-\int_t^T Z_s\,dW_s,\qquad 0\le t\le T, \tag{9}$$

*has a unique solution, that is, there exists a unique pair of $\mathscr{F}_t$-adapted processes $(Y_\cdot,Z_\cdot)\in S^2_{\mathscr{F}}(0,T;\mathbb{R}^m)\times L^2_{\mathscr{F}}(0,T;\mathbb{R}^{m\times d})$ satisfying equation (9).*

Lemma 3.3. *Assume $g_j(\omega,t,y,z):\Omega\times[0,T]\times\mathbb{R}\times\mathbb{R}^d\longrightarrow\mathbb{R}$ satisfies* (a) *and* (b), *$j=1,2$. Let $(Y^{(1)}_\cdot,Z^{(1)}_\cdot)$ and $(Y^{(2)}_\cdot,Z^{(2)}_\cdot)$ be respectively the solutions of BSDEs as follows:*

$$Y^{(j)}_t=\xi^{(j)}+\int_t^T g_j(s,Y^{(j)}_s,Z^{(j)}_s)\,ds-\int_t^T Z^{(j)}_s\,dW_s,\qquad 0\le t\le T,$$

*where $j=1,2$. If $\xi^{(1)}\ge\xi^{(2)}$ and $g_1(t,Y^{(1)}_t,Z^{(1)}_t)\ge g_2(t,Y^{(1)}_t,Z^{(1)}_t)$, a.e., a.s., then*

$$Y^{(1)}_t\ge Y^{(2)}_t,\qquad \text{a.e., a.s.}$$

*We also have strict comparison: under the above conditions,*

$$\begin{aligned}Y^{(1)}_0=Y^{(2)}_0\quad&\Longleftrightarrow\quad \xi^{(1)}=\xi^{(2)},\qquad \text{a.s.},\\ g_1(t,Y^{(1)}_t,Z^{(1)}_t)&=g_2(t,Y^{(1)}_t,Z^{(1)}_t),\qquad \text{a.e., a.s.}\end{aligned}$$

Lemma 3.4. *We make the same assumption as in Lemma 3.3. If $\xi^{(1)}\ge\xi^{(2)}$, $g_1(t,y,z)\ge g_2(t,y,z), t\in[0,T], y\in\mathbb{R}, z\in\mathbb{R}^d$, then*

$$Y^{(1)}_t\ge Y^{(2)}_t,\qquad \text{a.e., a.s.}$$

**4. Existence and uniqueness theorem.** We consider a new form of BSDEs as follows:

$$\begin{cases}-dY_t=f(t,Y_t,Z_t,Y_{t+\delta(t)},Z_{t+\zeta(t)})\,dt-Z_t\,dW_t, & t\in[0,T];\\ Y_t=\xi_t, & t\in[T,T+K];\\ Z_t=\eta_t, & t\in[T,T+K],\end{cases} \tag{10}$$

where $\delta(\cdot)$ and $\zeta(\cdot)$ are two $\mathbb{R}^+$-valued continuous functions defined on $[0,T]$ such that:

(i) There exists a constant $K\geq 0$ such that, for all $s\in[0,T]$,

$$s+\delta(s)\leq T+K;\qquad s+\zeta(s)\leq T+K.$$

(ii) There exists a constant $L\geq 0$ such that, for all $t\in[0,T]$ and for all nonnegative and integrable $g(\cdot)$,

$$\int_t^T g(s+\delta(s))\,ds\leq L\int_t^{T+K} g(s)\,ds;$$
$$\int_t^T g(s+\zeta(s))\,ds\leq L\int_t^{T+K} g(s)\,ds.$$

We call equation (10) the anticipated BSDE.

The setting of our problem is as follows: to find a pair of $\mathscr{F}_t$-adapted processes $(Y_\cdot,Z_\cdot)\in S^2_{\mathscr{F}}(0,T+K;\mathbb{R}^m)\times L^2_{\mathscr{F}}(0,T+K;\mathbb{R}^{m\times d})$ satisfying anticipated BSDE (10).

Assume that for all $s\in[0,T]$, $f(s,\omega,y,z,\xi,\eta):\Omega\times\mathbb{R}^m\times\mathbb{R}^{m\times d}\times L^2(\mathscr{F}_r;\mathbb{R}^m)\times L^2(\mathscr{F}_{r'};\mathbb{R}^{m\times d})\longrightarrow L^2(\mathscr{F}_s,\mathbb{R}^m)$, where $r,r'\in[s,T+K]$, and $f$ satisfies the following conditions:

(H1) There exists a constant $C>0$, such that for all $s\in[0,T], y,y'\in\mathbb{R}^m,\ z,z'\in\mathbb{R}^{m\times d},\xi_\cdot,\xi'_\cdot\in L^2_{\mathscr{F}}(s,T+K;\mathbb{R}^m),\eta_\cdot,\eta'_\cdot\in L^2_{\mathscr{F}}(s,T+K;\mathbb{R}^{m\times d}),r,$ $\bar{r}\in[s,T+K]$, we have

$$\begin{aligned}&|f(s,y,z,\xi_r,\eta_{\bar{r}})-f(s,y',z',\xi'_r,\eta'_{\bar{r}})|\\&\quad\leq C(|y-y'|+|z-z'|+E^{\mathscr{F}_s}[|\xi_r-\xi'_r|+|\eta_{\bar{r}}-\eta'_{\bar{r}}|]).\end{aligned}$$

(H2) $E[\int_0^T|f(s,0,0,0,0)|^2\,ds]<\infty$.

Remark 4.1. 1. Note that $f(s,\cdot,\cdot,\cdot,\cdot)$ is $\mathscr{F}_s$-measurable ensures the solution to the anticipated BSDE is $\mathscr{F}_s$-adapted.

2. We give examples of $\delta(s)$ and $f$. Both examples of $\delta(s)$ satisfy (i) and (ii). Example 1: Let $\delta(s)\equiv c$, where $c>0$ is a constant. Example 2: Let $s+\delta(s)$ be a monotone nonnegative function whose converse function has a continuous differential function. We give examples of functions that satisfy (H1) and (H2): Let $g$ satisfy (a) and (b) and let $\delta,\zeta$ be two positive constants. For each $t\in[0,T]$ and $(\xi_\cdot,\eta_\cdot)\in L^2_{\mathscr{F}}(t,T+(\delta\vee\zeta);\mathbb{R}^m\times\mathbb{R}^{m\times d})$, define $f_1,f_2$ such that

$$f_1(t,\xi_{t+\delta},\eta_{t+\zeta})=g(t,E^{\mathscr{F}_t}[\xi_{t+\delta}],E^{\mathscr{F}_t}[\eta_{t+\zeta}]),$$
$$f_2(t,\xi_{t+\delta},\eta_{t+\zeta})=g(t,\mathcal{E}_{t,t+\delta}[\xi_{t+\delta}],\mathcal{E}_{t,t+\zeta}[\eta_{t+\zeta}]),$$

where $\mathcal{E}_{s,t}[\cdot]:L^2(\mathscr{F}_t)\longrightarrow L^2(\mathscr{F}_s), 0\leq s\leq t\leq T+K$, is a $\mathscr{F}_t$-consistent nonlinear evaluation (see Peng [13]). Then $f_1,f_2$ satisfy (i) and (ii).

The following is the main result of this section: an existence and uniqueness theorem for adapted solutions for anticipated BSDEs.

THEOREM 4.2. *Suppose that $f$ satisfies* (H1) *and* (H2), *and $\delta, \zeta$ satisfy* (i) *and* (ii). *Then for any given terminal conditions $\xi_\bullet \in S^2_{\mathscr{F}}(T, T+K; \mathbb{R}^m)$ and $\eta_\bullet \in L^2_{\mathscr{F}}(T, T+K; \mathbb{R}^{m\times d})$, the anticipated BSDE (10) has a unique solution, that is, there exists a unique pair of $\mathscr{F}_t$-adapted processes $(Y_\cdot, Z_\cdot) \in S^2_{\mathscr{F}}(0, T+K; \mathbb{R}^m) \times L^2_{\mathscr{F}}(0, T+K; \mathbb{R}^{m\times d})$ satisfying (10).*

PROOF. We fix $\beta = 12C^2(2L+1) + 2$, where $C$ is the Lipschitz constant of $f$ given in (H1), and introduce a norm in the Banach space $L^2_{\mathscr{F}}(0, T+K; \mathbb{R}^m)$:

$$\|\nu(\cdot)\|_\beta = \left(E\left[\int_0^{T+K} |\nu_s|^2 e^{\beta s}\, ds\right]\right)^{1/2}.$$

Clearly, it is equivalent to the original norm of $L^2_{\mathscr{F}}(0, T+K; \mathbb{R}^m)$. But it is more convenient to use this norm to construct a contraction mapping that allows us to apply the Fixed Point Theorem. Set

$$\begin{cases} Y_t = \xi_T + \displaystyle\int_t^T f(s, y_s, z_s, y_{s+\delta(s)}, z_{s+\zeta(s)})\, ds - \int_t^T Z_s\, dW_s, \\ \qquad\qquad t \in [0, T]; \\ Y_t = \xi_t, \qquad t \in [T, T+K]; \\ Z_t = \eta_t, \qquad t \in [T, T+K]. \end{cases}$$

Define a mapping $h : L^2_{\mathscr{F}}(0, T+K; \mathbb{R}^m \times \mathbb{R}^{m\times d}) \longrightarrow L^2_{\mathscr{F}}(0, T+K; \mathbb{R}^m \times \mathbb{R}^{m\times d})$ such that $h[(y_\bullet, z_\bullet)] = (Y_\bullet, Z_\bullet)$. Now we prove that $h$ is a contraction mapping under the norm $\|\cdot\|_\beta$. For two arbitrary elements $(y_\bullet, z_\bullet)$ and $(y'_\bullet, z'_\bullet)$ in $L^2_{\mathscr{F}}(0, T+K; \mathbb{R}^m \times \mathbb{R}^{m\times d})$, set $(Y_\bullet, Z_\bullet) = h[(y_\bullet, z_\bullet)]$ and $(Y'_\bullet, Z'_\bullet) = h[(y'_\bullet, z'_\bullet)]$. Denote their differences by

$$(\hat{y}_\bullet, \hat{z}_\bullet) = ((y - y')_\bullet, (z - z')_\bullet), \qquad (\hat{Y}_\bullet, \hat{Z}_\bullet) = ((Y - Y')_\bullet, (Z - Z')_\bullet).$$

By basic estimate (7), we have

$$\begin{aligned} &E\left[\int_0^T \left(\frac{\beta}{2}|\hat{Y}_s|^2 + |\hat{Z}_s|^2\right) e^{\beta s}\, ds\right] \\ &\qquad \le \frac{2}{\beta} E\left[\int_0^T |f(s, y_s, z_s, y_{s+\delta(s)}, z_{s+\zeta(s)}) \right. \\ &\qquad\qquad\qquad \left. - f(s, y'_s, z'_s, y'_{s+\delta(s)}, z'_{s+\zeta(s)})|^2 e^{\beta s}\, ds\right]. \end{aligned}$$

Since $\delta(s)$ and $\zeta(s)$ satisfy (ii) and $f$ satisfies (H1), by the Fubini Theorem, we have

$$\begin{aligned}
&E\biggl[\int_0^T \biggl(\frac{\beta}{2}|\hat{Y}_s|^2 + |\hat{Z}_s|^2\biggr) e^{\beta s}\, ds\biggr] \\
&\quad \le \frac{2C^2}{\beta} E\biggl[\int_0^T (|\hat{y}_s| + |\hat{z}_s| + E^{\mathscr{F}_s}[|\hat{y}_{s+\delta(s)}| + |\hat{z}_{s+\zeta(s)}|])^2 e^{\beta s}\, ds\biggr] \\
&\quad \le \frac{6C^2}{\beta} E\biggl[\int_0^T (|\hat{y}_s|^2 + |\hat{z}_s|^2 + 2|\hat{y}_{s+\delta(s)}|^2 + 2|\hat{z}_{s+\zeta(s)}|^2) e^{\beta s}\, ds\biggr] \\
&\quad \le \frac{6C^2(2L+1)}{\beta} E\biggl[\int_0^{T+K} (|\hat{y}_s|^2 + |\hat{z}_s|^2) e^{\beta s}\, ds\biggr].
\end{aligned}$$

Because $\beta = 12C^2(2L+1) + 2$, then

$$\begin{aligned}
&E\biggl[\int_0^{T+K} (|\hat{Y}_s|^2 + |\hat{Z}_s|^2) e^{\beta s}\, ds\biggr] \\
&\qquad \le \frac{1}{2} E\biggl[\int_0^{T+K} (|\hat{y}_s|^2 + |\hat{z}_s|^2) e^{\beta s}\, ds\biggr],
\end{aligned}$$

or

$$\|(\hat{Y}_\bullet, \hat{Z}_\bullet)\|_\beta \le \tfrac{1}{\sqrt{2}} \|(\hat{y}_\bullet, \hat{z}_\bullet)\|_\beta.$$

Consequently, $h$ is a strict contraction mapping of $L^2_{\mathscr{F}}(0, T+K; \mathbb{R}^m \times \mathbb{R}^{m\times d})$. It follows by the Fixed Point Theorem that (10) has a unique solution $(Y_\bullet, Z_\bullet) \in L^2_{\mathscr{F}}(0, T+K; \mathbb{R}^m \times \mathbb{R}^{m\times d})$. Since $f$ satisfies (H1) and (H2) and since $\delta, \zeta$ satisfy (i) and (ii), we have $f(\cdot, Y_\bullet, Z_\bullet, Y_{\bullet+\delta(\bullet)}, Z_{\bullet+\zeta(\bullet)}) \in L^2_{\mathscr{F}}(0, T; \mathbb{R}^m)$. Thus, by Lemma 3.1, we obtain $Y_\bullet \in S^2_{\mathscr{F}}(0, T+K; \mathbb{R}^m)$. $\square$

The following example shows that a simple case of the anticipated BSDE (10) has a solution.

Example 4.3. Consider the following anticipated BSDE:

$$\begin{cases} Y_t = TW_T - \displaystyle\int_t^T \frac{1}{s+\delta} E^{\mathscr{F}_s}[Y_{s+\delta}]\, ds - \int_t^T Z_s\, dW_s, & t \in [0, T]; \\ Y_t = tW_t, & t \in [T, T+\delta], \end{cases}$$

where $\delta \ge 0$ is a given constant. Then $(tW_t, t)_{t\in[0,T+\delta]}$ is its solution.

The following proposition is an estimate of the solution of the anticipated BSDE (10).

Proposition 4.4. *Assume that $f$ satisfies* (H1) *and* (H2), *and also $\delta$ and $\zeta$ satisfy* (i) *and* (ii). *Then there exists a positive constant $C_0$ that only* *depends on $C$ in* (H1), *$L$ in* (ii), *and $T$ such that for each $\xi_\bullet \in S^2_{\mathscr{F}}(T,T+K;\mathbb{R}^m)$ and each $\eta_\bullet \in L^2_{\mathscr{F}}(T,T+K;\mathbb{R}^{m\times d})$, the solution $(Y_\bullet,Z_\bullet)$ of the anticipated BSDE (10) satisfies*

$$\begin{aligned} &E^{\mathscr{F}_t}\left[\sup_{t\le s\le T}|Y_s|^2+\int_t^T|Z_s|^2\,ds\right] \\ &\quad\le C_0E^{\mathscr{F}_t}\left[|\xi_T|^2+\int_T^{T+K}(|\xi_s|^2+|\eta_s|^2)\,ds+\left(\int_t^T|f(s,0,0,0,0)|\,ds\right)^2\right], \end{aligned}\tag{11}$$

*for each $t\in[0,T]$.*

Proof. For $s\in[0,T]$, applying Itô's formula to $e^{\beta s}|Y_s|^2$, we obtain

$$\begin{aligned} &e^{\beta s}|Y_s|^2+\int_s^T e^{\beta r}(\beta|Y_r|^2+|Z_r|^2)\,dr \\ &\quad= e^{\beta T}|\xi_T|^2-2\int_s^T e^{\beta r}(Y_r,Z_r\,dW_r) \\ &\qquad+2\int_s^T e^{\beta r}(f(r,Y_r,Z_r,Y_{r+\delta(r)},Z_{r+\zeta(r)}),Y_r)\,dr. \end{aligned}$$

Since

$$\begin{aligned} &2(f(r,Y_r,Z_r,Y_{r+\delta(r)},Z_{r+\zeta(r)}),Y_r) \\ &\quad=2(f(r,Y_r,Z_r,Y_{r+\delta(r)},Z_{r+\zeta(r)})-f(r,Y_r,Z_r,Y_{r+\delta(r)},0),Y_r) \\ &\qquad+2(f(r,Y_r,Z_r,Y_{r+\delta(r)},0)-f(r,Y_r,Z_r,0,0),Y_r) \\ &\qquad+2(f(r,Y_r,Z_r,0,0)-f(r,Y_r,0,0,0),Y_r) \\ &\qquad+2(f(r,Y_r,0,0,0)-f(r,0,0,0,0),Y_r)+2(f(r,0,0,0,0),Y_r) \\ &\quad\le 2CE^{\mathscr{F}_r}[|Z_{r+\zeta(r)}|]|Y_r|+2CE^{\mathscr{F}_r}[|Y_{r+\delta(r)}|]|Y_r|+2C|Y_r||Z_r| \\ &\qquad+2C|Y_r|^2+2(f(r,0,0,0,0),Y_r) \\ &\quad\le(3LC^2+4LTC^2+3C^2+2C)|Y_r|^2+\frac{1}{3L}E^{\mathscr{F}_r}[|Z_{r+\zeta(r)}|^2] \\ &\qquad+\frac{1}{4LT}E^{\mathscr{F}_r}[|Y_{r+\delta(r)}|^2]+\frac{1}{3}|Z_r|^2+2(f(r,0,0,0,0),Y_r), \end{aligned}$$

we get, for $s\in[0,T]$,

$$\begin{aligned} &e^{\beta s}|Y_s|^2+\int_s^T e^{\beta r}\left[(\beta-3LC^2-4LTC^2-3C^2-2C)|Y_r|^2+\frac{2}{3}|Z_r|^2\right]dr \\ &\quad\le e^{\beta T}|\xi_T|^2+2\int_s^T e^{\beta r}(f(r,0,0,0,0),Y_r)\,dr-2\int_s^T e^{\beta r}(Y_r,Z_r\,dW_r) \end{aligned}\tag{12}$$

$$+\frac{1}{3L}\int_s^T e^{\beta r}E^{\mathscr{F}_r}[|Z_{r+\zeta(r)}|^2]\,dr+\frac{1}{4LT}\int_s^T e^{\beta r}E^{\mathscr{F}_r}[|Y_{r+\delta(r)}|^2]\,dr.$$

Taking conditional expectations under $\mathscr{F}_s$ on both sides of (12), we have

$$\begin{aligned}
&e^{\beta s}|Y_s|^2+E^{\mathscr{F}_s}\biggl[\int_s^T e^{\beta r}\biggl[(\beta-3LC^2-4LTC^2-3C^2-2C)|Y_r|^2\\
&\qquad\qquad\qquad\qquad\qquad\qquad+\frac{2}{3}|Z_r|^2\biggr]\,dr\biggr]\\
&\quad\le E^{\mathscr{F}_s}\biggl[e^{\beta T}|\xi_T|^2+2\int_s^T e^{\beta r}(f(r,0,0,0,0),Y_r)\,dr\biggr]\\
&\qquad+\frac{1}{4LT}E^{\mathscr{F}_s}\biggl[\int_s^T e^{\beta r}E^{\mathscr{F}_r}[|Y_{r+\delta(r)}|^2]\,dr\biggr]\\
&\qquad+\frac{1}{3L}E^{\mathscr{F}_s}\biggl[\int_s^T e^{\beta r}E^{\mathscr{F}_r}[|Z_{r+\zeta(r)}|^2]\,dr\biggr]\\
&\quad\le E^{\mathscr{F}_s}\biggl[e^{\beta T}|\xi_T|^2+2\int_s^T e^{\beta r}(f(r,0,0,0,0),Y_r)\,dr\\
&\qquad\qquad\qquad+\frac{1}{4T}\int_s^{T+K}e^{\beta r}|Y_r|^2\,dr\biggr]\\
&\qquad+\frac{1}{3}E^{\mathscr{F}_s}\biggl[\int_s^{T+K}e^{\beta r}|Z_r|^2\,dr\biggr].
\end{aligned}$$

Set $\beta=3LC^2+4LTC^2+3C^2+2C+\dfrac{1}{4T}$, then

$$\begin{aligned}
&E^{\mathscr{F}_s}\biggl[\int_s^T e^{\beta r}|Z_r|^2\,dr\biggr]\\
&\quad\le E^{\mathscr{F}_s}\biggl[3e^{\beta T}|\xi_T|^2+6\int_s^T e^{\beta r}(f(r,0,0,0,0),Y_r)\,dr\biggr]\\
&\qquad+E^{\mathscr{F}_s}\biggl[\int_T^{T+K}e^{\beta r}\biggl(|\eta_r|^2+\frac{3}{4T}|\xi_r|^2\biggr)\,dr\biggr].
\end{aligned}\tag{13}$$

Since for $t\le s\le T$,

$$\begin{aligned}
\biggl|\int_s^T e^{\beta r}(Y_r,Z_r\,dW_r)\biggr|&=\biggl|\int_t^T e^{\beta r}(Y_r,Z_r\,dW_r)-\int_t^s e^{\beta r}(Y_r,Z_r\,dW_r)\biggr|\\
&\le\biggl|\int_t^T e^{\beta r}(Y_r,Z_r\,dW_r)\biggr|+\biggl|\int_t^s e^{\beta r}(Y_r,Z_r\,dW_r)\biggr|,
\end{aligned}$$

by the Burkholder–Davis–Gundy inequality, we have

$$E^{\mathscr{F}_t}\biggl[\sup_{t\le s\le T}\biggl|\int_s^T e^{\beta r}(Y_r,Z_r\,dW_r)\biggr|\biggr]$$

$$
\begin{aligned}
&\le 2E^{\mathscr{F}_t}\left[\sup_{t\le s\le T}\left|\int_t^s e^{\beta r}(Y_r, Z_r\,dW_r)\right|\right]\\
&\le 6E^{\mathscr{F}_t}\left[\left(\int_t^T e^{2\beta r}|Y_r|^2|Z_r|^2\,dr\right)^{1/2}\right] \qquad (14)\\
&\le 6E^{\mathscr{F}_t}\left[\left(\sup_{t\le r\le T} e^{1/2\beta r}|Y_r|\right)\left(\int_t^T e^{\beta r}|Z_r|^2\,dr\right)^{1/2}\right]\\
&\le \frac{1}{4}E^{\mathscr{F}_t}\left[\sup_{t\le r\le T} e^{\beta r}|Y_r|^2\right] + 36E^{\mathscr{F}_t}\left[\int_t^T e^{\beta r}|Z_r|^2\,dr\right].
\end{aligned}
$$

From estimates (12) and (14) we have

$$
\begin{aligned}
&E^{\mathscr{F}_t}\left[\sup_{t\le s\le T} e^{\beta s}|Y_s|^2\right]\\
&\quad\le E^{\mathscr{F}_t}\Big[e^{\beta T}|\xi_T|^2 + 2\int_t^T e^{\beta r}|f(r,0,0,0,0)||Y_r|\,dr\\
&\qquad\qquad + 2\sup_{t\le s\le T}\left|\int_s^T e^{\beta r}(Y_r, Z_r\,dW_r)\right|\Big]\\
&\qquad + E^{\mathscr{F}_t}\Big[\frac{1}{3L}\int_t^T e^{\beta r}E^{\mathscr{F}_r}[|Z_{r+\zeta(r)}|^2]\,dr\\
&\qquad\qquad + \frac{1}{4LT}\int_t^T e^{\beta r}E^{\mathscr{F}_r}[|Y_{r+\delta(r)}|^2]\,dr\Big]\\
&\quad\le E^{\mathscr{F}_t}[e^{\beta T}|\xi_T|^2] + \frac{1}{2}E^{\mathscr{F}_t}\left[\sup_{t\le r\le T} e^{\beta r}|Y_r|^2\right]\\
&\qquad + 72E^{\mathscr{F}_t}\left[\int_t^T e^{\beta r}|Z_r|^2\,dr\right]\\
&\qquad + E^{\mathscr{F}_t}\left[\frac{1}{3L}\int_t^T e^{\beta r}|Z_{r+\zeta(r)}|^2\,dr + \frac{1}{4LT}\int_t^T e^{\beta r}|Y_{r+\delta(r)}|^2\,dr\right]\\
&\qquad + 2E^{\mathscr{F}_t}\left[\int_t^T e^{\beta r}|f(r,0,0,0,0)||Y_r|\,dr\right]\\
&\quad\le E^{\mathscr{F}_t}[e^{\beta T}|\xi_T|^2] + \frac{1}{2}E^{\mathscr{F}_t}\left[\sup_{t\le r\le T} e^{\beta r}|Y_r|^2\right] + 72E^{\mathscr{F}_t}\left[\int_t^T e^{\beta r}|Z_r|^2\,dr\right]\\
&\qquad + E^{\mathscr{F}_t}\left[\frac{1}{3}\int_t^{T+K} e^{\beta r}|Z_r|^2\,dr + \frac{1}{4T}\int_t^{T+K} e^{\beta r}|Y_r|^2\,dr\right]\\
&\qquad + 2E^{\mathscr{F}_t}\left[\int_t^T e^{\beta r}|f(r,0,0,0,0)||Y_r|\,dr\right]
\end{aligned}
$$

$$
\begin{aligned}
&\leq E^{\mathscr{F}_t}[e^{\beta T}|\xi_T|^2] + \frac{3}{4}E^{\mathscr{F}_t}\left[\sup_{t\leq r\leq T} e^{\beta r}|Y_r|^2\right] \\
&\quad + \left(72+\frac{1}{3}\right)E^{\mathscr{F}_t}\left[\int_t^T e^{\beta r}|Z_r|^2\,dr\right] \\
&\quad + E^{\mathscr{F}_t}\left[\int_T^{T+K} e^{\beta r}\left(\frac{1}{3}|\eta_r|^2 + \frac{1}{4T}|\xi_r|^2\right)dr\right] \\
&\quad + 2E^{\mathscr{F}_t}\left[\int_t^T e^{\beta r}|f(r,0,0,0,0)||Y_r|\,dr\right].
\end{aligned}
$$

Denote by $C_0 > 0$ a constant that depends only on $T, L$ and $C$, which we allow to change from line to line. From the estimate above and estimate (13),

$$
\begin{aligned}
&\frac{1}{4}E^{\mathscr{F}_t}\left[\sup_{t\leq s\leq T} e^{\beta s}|Y_s|^2\right] \\
&\quad\leq C_0 E^{\mathscr{F}_t}\left[e^{\beta T}|\xi_T|^2 + \int_T^{T+K} e^{\beta r}(|\eta_r|^2 + |\xi_r|^2)\,dr\right] \\
&\qquad + C_0 E^{\mathscr{F}_t}\left[\left(\sup_{t\leq r\leq T} e^{1/2\beta r}|Y_r|\right)\left(\int_t^T e^{1/2\beta r}|f(r,0,0,0,0)|\,dr\right)\right] \\
&\quad\leq C_0 E^{\mathscr{F}_t}\left[e^{\beta T}|\xi_T|^2 + \int_T^{T+K} e^{\beta r}(|\eta_r|^2 + |\xi_r|^2)\,dr\right] \\
&\qquad + \frac{1}{8}E^{\mathscr{F}_t}\left[\sup_{t\leq r\leq T} e^{\beta r}|Y_r|^2\right] + 2C_0^2 E^{\mathscr{F}_t}\left[\left(\int_t^T e^{1/2\beta r}|f(r,0,0,0,0)|\,dr\right)^2\right] \\
&\quad\leq C_0 E^{\mathscr{F}_t}\left[|\xi_T|^2 + \int_T^{T+K}(|\eta_r|^2 + |\xi_r|^2)\,dr\right] \\
&\qquad + \frac{1}{8}E^{\mathscr{F}_t}\left[\sup_{t\leq r\leq T} e^{\beta r}|Y_r|^2\right] + 2C_0^2 E^{\mathscr{F}_t}\left[\left(\int_t^T |f(r,0,0,0,0)|\,dr\right)^2\right].
\end{aligned}
$$

Then

$$
\begin{aligned}
&E^{\mathscr{F}_t}\left[\sup_{t\leq s\leq T}|Y_s|^2\right] + E^{\mathscr{F}_t}\left[\int_t^T |Z_s|^2\,ds\right] \\
&\quad\leq C_0 E^{\mathscr{F}_t}\left[|\xi_T|^2 + \int_T^{T+K}(|\xi_s|^2 + |\eta_s|^2)\,ds + \left(\int_t^T |f(s,0,0,0,0)|\,ds\right)^2\right].
\end{aligned}
$$

□

The following proposition shows the importance of the effect of anticipated time on the solution to anticipated BSDEs.

PROPOSITION 4.5. *Let $(Y_\cdot^{(1)}, Z_\cdot^{(1)})$ and $(Y_\cdot^{(2)}, Z_\cdot^{(2)})$ be respectively solutions of the following two anticipated BSDEs:*

$$\begin{cases} -dY_t^{(j)} = f(t, Y_t^{(j)}, Z_t^{(j)}, Y_{t+\delta_j(t)}^{(j)})\,dt - Z_t^{(j)}\,dW_t, & t \in [0,T]; \\ Y_t^{(j)} = \xi_t, & t \in [T, T+K], \end{cases}$$

*where $j = 1,2$. Assume $\xi_\cdot \in S_{\mathscr{F}}^2(T, T+K; \mathbb{R}^m)$, $\delta_1$ and $\delta_2$ satisfy* (i) *and* (ii), *$f$ satisfies* (H2), *and there exists a constant $\bar{C} > 0$, such that for all $s \in [0,T], y, y' \in \mathbb{R}^m, z, z' \in \mathbb{R}^{m\times d}, \theta, \theta' \in L_{\mathscr{F}}^2(s, T+K; \mathbb{R}^m)$ and $r \in [s, T+K]$,*

$$|f(s,y,z,\theta_r) - f(s,y',z',\theta_r')| \le \bar{C}(|y-y'| + |z-z'| + |E^{\mathscr{F}_s}[\theta_r - \theta_r']|).$$

*If for any $t \in [0,T]$, $\delta_1(t) \le \delta_2(t)$, then there exists a constant $\tilde{M} > 0$ only depending on $\bar{C}$, $L$ and $T$ such that*

$$|Y_t^{(1)} - Y_t^{(2)}|^2 \le \tilde{M}\int_t^T (\delta_2(s) - \delta_1(s))\,ds \times E^{\mathscr{F}_t}\left[|\xi_T|^2 + \int_T^{T+K} |\xi_s|^2\,ds + \int_t^T |f(s,0,0,0)|^2\,ds\right].$$

PROOF. Setting $y_\cdot = Y_\cdot^{(1)} - Y_\cdot^{(2)}$, $z_\cdot = Z_\cdot^{(1)} - Z_\cdot^{(2)}$, then by estimate (6), we obtain, for all $t \in [0,T]$,

$$\begin{aligned}
&|y_t|^2 + E^{\mathscr{F}_t}\left[\int_t^T \left(\frac{\beta}{2}|y_s|^2 + |z_s|^2\right)e^{\beta(s-t)}\,ds\right] \\
&\quad\le \frac{2}{\beta}E^{\mathscr{F}_t}\left[\int_t^T |f(s, Y_s^{(1)}, Z_s^{(1)}, Y_{s+\delta_1(s)}^{(1)}) - f(s, Y_s^{(2)}, Z_s^{(2)}, Y_{s+\delta_2(s)}^{(2)})|^2 \times e^{\beta(s-t)}\,ds\right] \\
&\quad\le \frac{2\bar{C}^2}{\beta}E^{\mathscr{F}_t}\left[\int_t^T (|y_s| + |z_s| + |E^{\mathscr{F}_s}[Y_{s+\delta_1(s)}^{(1)} - Y_{s+\delta_2(s)}^{(2)}]|)^2 e^{\beta(s-t)}\,ds\right] \\
&\quad\le \frac{6\bar{C}^2}{\beta}E^{\mathscr{F}_t}\left[\int_t^T (|y_s|^2 + |z_s|^2 + |E^{\mathscr{F}_s}[y_{s+\delta_1(s)}] + E^{\mathscr{F}_s}[Y_{s+\delta_1(s)}^{(2)} - Y_{s+\delta_2(s)}^{(2)}]|^2)e^{\beta(s-t)}\,ds\right] \\
&\quad\le \frac{6\bar{C}^2 + 12\bar{C}^2 L}{\beta}E^{\mathscr{F}_t}\left[\int_t^T |y_s|^2 e^{\beta(s-t)}\,ds\right] + \frac{6\bar{C}^2}{\beta}E^{\mathscr{F}_t}\left[\int_t^T |z_s|^2 e^{\beta(s-t)}\,ds\right] \\
&\qquad + \frac{12\bar{C}^2}{\beta}E^{\mathscr{F}_t}\left[\int_t^T |E^{\mathscr{F}_s}[Y_{s+\delta_1(s)}^{(2)} - Y_{s+\delta_2(s)}^{(2)}]|^2 e^{\beta(s-t)}\,ds\right].
\end{aligned}$$

But

$$E^{\mathscr{F}_s}[Y^{(2)}_{s+\delta_1(s)} - Y^{(2)}_{s+\delta_2(s)}] = E^{\mathscr{F}_s}\left[\int_{s+\delta_1(s)}^{s+\delta_2(s)} f(r, Y_r^{(2)}, Z_r^{(2)}, Y^{(2)}_{r+\delta_2(r)})\,dr\right],$$

and set $\beta = 6\bar{C}^2$, hence,

$$\begin{aligned}
|y_t|^2 &\le (1+2L)E^{\mathscr{F}_t}\left[\int_t^T |y_s|^2 e^{\beta(s-t)}\,ds\right] \\
&\quad + 2E^{\mathscr{F}_t}\left[\int_t^T E^{\mathscr{F}_s}\left[\left|\int_{s+\delta_1(s)}^{s+\delta_2(s)} f(r, Y_r^{(2)}, Z_r^{(2)}, Y^{(2)}_{r+\delta_2(r)})\,dr\right|^2\right] e^{\beta(s-t)}\,ds\right] \\
&\le (1+2L)e^{\beta(T-t)}E^{\mathscr{F}_t}\left[\int_t^T |y_s|^2\,ds\right] \\
&\quad + 2E^{\mathscr{F}_t}\left[\int_t^T (\delta_2(s)-\delta_1(s))\int_{s+\delta_1(s)}^{s+\delta_2(s)} |f(r, Y_r^{(2)}, Z_r^{(2)}, Y^{(2)}_{r+\delta_2(r)})|^2\,dr\right. \\
&\qquad\qquad\qquad\qquad\qquad\qquad\qquad\qquad \left.\times e^{\beta(s-t)}\,ds\right] \\
&\le (1+2L)e^{\beta(T-t)}E^{\mathscr{F}_t}\left[\int_t^T |y_s|^2\,ds\right] \\
&\quad + 8E^{\mathscr{F}_t}\left[\int_t^T (\delta_2(s)-\delta_1(s))e^{\beta(s-t)}\,ds\right. \\
&\qquad\qquad \times\left\{\int_t^T (\bar{C}^2|Y_r^{(2)}|^2 + \bar{C}^2|Z_r^{(2)}|^2 + \bar{C}^2|Y^{(2)}_{r+\delta_2(r)}|^2\right. \\
&\qquad\qquad\qquad\qquad\qquad\qquad \left.\left. + |f(r,0,0,0)|^2)\,dr\right\}\right] \\
&\le (1+2L)e^{\beta(T-t)}E^{\mathscr{F}_t}\left[\int_t^T |y_s|^2\,ds\right] \\
&\quad + 8\int_t^T (\delta_2(s)-\delta_1(s))e^{\beta(s-t)}\,ds \\
&\quad \times E^{\mathscr{F}_t}\left[\int_t^T ((1+L)\bar{C}^2|Y_r^{(2)}|^2 + \bar{C}^2|Z_r^{(2)}|^2 + |f(r,0,0,0)|^2)\,dr\right. \\
&\qquad\qquad\qquad\qquad\qquad \left. + \int_T^{T+K} L\bar{C}^2|\xi_r|^2\,dr\right].
\end{aligned}$$

From estimate (11), we can find a constant $\bar{M} > 0$ depending only on $\bar{C}$, $L$ and $T$ such that

$$|y_t|^2 \le \bar{M}E^{\mathscr{F}_t}\left[\int_t^T |y_s|^2\,ds\right]$$

$$
\begin{aligned}
&+ \bar{M} \int_t^T (\delta_2(s) - \delta_1(s))\, ds \\
&\quad \times E^{\mathscr{F}_t}\left[|\xi_T|^2 + \int_T^{T+K} |\xi_s|^2\, ds + \int_t^T |f(s,0,0,0)|^2\, ds\right].
\end{aligned}
$$

Thus, by Gronwall's inequality,

$$
\begin{aligned}
|y_t|^2 \le{}& \bar{M} \int_t^T (\delta_2(s) - \delta_1(s))\, ds \\
&\times E^{\mathscr{F}_t}\left[|\xi_T|^2 + \int_T^{T+K} |\xi_s|^2\, ds + \int_t^T |f(s,0,0,0)|^2\, ds\right] e^{\bar{M}t}.
\end{aligned}
$$

Fix $\tilde{M} = \bar{M} e^{\bar{M}T}$, therefore,

$$
\begin{aligned}
|y_t|^2 \le{}& \tilde{M} \int_t^T (\delta_2(s) - \delta_1(s))\, ds \\
&\times E^{\mathscr{F}_t}\left[|\xi_T|^2 + \int_T^{T+K} |\xi_s|^2\, ds + \int_t^T |f(s,0,0,0)|^2\, ds\right]. \qquad \square
\end{aligned}
$$

**5. Comparison theorem for 1-dimensional anticipated BSDEs.** Lemma 3.3 is a typical version of a comparison theorem. It is a fundamentally important result in BSDE theory. Some further developments in this direction are Cao and Yan [3], Lin [10], Liu and Ren [11], Zhang [16] and Situ [15], without mentioning many other widely circulated papers listed in [13]. Recently Hu and Peng [8] gave a comparison theorem for multidimensional BSDEs. Comparison theorems for BSDEs have received a lot of attention because of their importance. For example, the punishment method in reflected BSDEs is based on a comparison theorem (see [4, 6, 9] and [14]). Moreover, research on properties of g-expectations (see Peng [13]) and the proof of a monotonic limit theorem for BSDEs (see Peng [12]) both depend on comparison theorems.

It is well known that 1-dimensional BSDEs have comparison theorems (see Lemmas 3.3 and 3.4) when their generators satisfy the conditions of existence and uniqueness theorems for BSDEs. It is very important to notice that the conditions on $f$ needed for the comparison theorem for anticipated BSDEs are stronger than those needed for the existence and uniqueness theorem. Using the comparison theorem for anticipated BSDEs, we will solve a stochastic control problem in Section 6.

Let $(Y_\cdot^{(1)}, Z_\cdot^{(1)}), (Y_\cdot^{(2)}, Z_\cdot^{(2)})$ be respectively solutions of the following two 1-dimensional anticipated BSDEs:

$$
\begin{cases}
-dY_t^{(j)} = f_j(t, Y_t^{(j)}, Z_t^{(j)}, Y_{t+\delta(t)}^{(j)})\, dt - Z_t^{(j)}\, dW_t, & 0 \le t \le T; \\
Y_t^{(j)} = \xi_t^{(j)}, & T \le t \le T+K,
\end{cases}
$$

where $j = 1, 2$.

THEOREM 5.1. *Assume that* $f_1, f_2$ *satisfies* (H1) *and* (H2), $\xi_\cdot^{(1)}, \xi_\cdot^{(2)} \in S^2_{\mathscr{F}}(T, T+K)$, $\delta$ *satisfies* (i), (ii), *and for all* $t \in [0,T], y \in \mathbb{R}, z \in \mathbb{R}^d$, $f_2(t, y, z, \cdot)$ *is increasing, that is,* $f_2(t,y,z,\theta_r) \geq f_2(t,y,z,\theta'_r)$, *if* $\theta_r \geq \theta'_r$, $\theta, \theta' \in L^2_{\mathscr{F}}(t, T+K), r \in [t, T+K]$. *If* $\xi_s^{(1)} \geq \xi_s^{(2)}, s \in [T, T+K]$ *and* $f_1(t,y,z,\theta_r) \geq f_2(t,y,z,\theta_r), t \in [0,T], y \in \mathbb{R}, z \in \mathbb{R}^d, \theta \in L^2_{\mathscr{F}}(t, T+K), r \in [t, T+K]$, *then*

$$Y_t^{(1)} \geq Y_t^{(2)}, \qquad a.e., \quad a.s.$$

PROOF. Set

$$\begin{cases} Y_t^{(3)} = \xi_T^{(2)} + \displaystyle\int_t^T f_2(s, Y_s^{(3)}, Z_s^{(3)}, Y_{s+\delta(s)}^{(1)})\, ds - \int_t^T Z_s^{(3)}\, dW_s, \\ \qquad\qquad\qquad t \in [0,T]; \\ Y_t^{(3)} = \xi_t^{(2)}, \qquad t \in [T, T+K]. \end{cases}$$

By Lemma 3.2, we know there exists a unique pair of $\mathscr{F}_t$-adapted processes $(Y_\cdot^{(3)}, Z_\cdot^{(3)}) \in S^2_{\mathscr{F}}(0,T) \times L^2_{\mathscr{F}}(0,T;\mathbb{R}^d)$ that satisfies the above BSDE. Since $f_1(s,y,z,Y_{s+\delta(s)}^{(1)}) \geq f_2(s,y,z,Y_{s+\delta(s)}^{(1)})$, $s \in [0,T], y \in \mathbb{R}, z \in \mathbb{R}^d$, by Lemma 3.4, we obtain

$$Y_t^{(1)} \geq Y_t^{(3)}, \qquad \text{a.e.,} \quad \text{a.s.}$$

Set

$$\begin{cases} Y_t^{(4)} = \xi_T^{(2)} + \displaystyle\int_t^T f_2(s, Y_s^{(4)}, Z_s^{(4)}, Y_{s+\delta(s)}^{(3)})\, ds - \int_t^T Z_s^{(4)}\, dW_s, \\ \qquad\qquad\qquad t \in [0,T]; \\ Y_t^{(4)} = \xi_t^{(2)}, \qquad t \in [T, T+K]. \end{cases}$$

Since for all $t \in [0,T], y \in \mathbb{R}, z \in \mathbb{R}^d$, $f_2(t,y,z,\cdot)$ is increasing and $Y_t^{(1)} \geq Y_t^{(3)}$, a.e., a.s., by Lemma 3.4, we know

$$Y_t^{(3)} \geq Y_t^{(4)}, \qquad \text{a.e.,} \quad \text{a.s.}$$

For $n = 5, 6, \ldots,$ we consider the following classical BSDE:

$$\begin{cases} Y_t^{(n)} = \xi_T^{(2)} + \displaystyle\int_t^T f_2(s, Y_s^{(n)}, Z_s^{(n)}, Y_{s+\delta(s)}^{(n-1)})\, ds - \int_t^T Z_s^{(n)}\, dW_s, \\ \qquad\qquad\qquad t \in [0,T]; \\ Y_t^{(n)} = \xi_t^{(2)}, \qquad t \in [T, T+K]. \end{cases}$$

Similarly, we have $Y_t^{(4)} \geq Y_t^{(5)} \geq \cdots \geq Y_t^{(n)} \geq \cdots$, a.e., a.s. We use $\|\nu(\cdot)\|_\beta$ in the proof of Theorem 4.2 as the norm in the Banach space $L^2_{\mathscr{F}}(0, T+K;\mathbb{R}) \times L^2_{\mathscr{F}}(0,T;\mathbb{R}^d)$. Set $\hat{Y}_s^{(n)} = Y_s^{(n)} - Y_s^{(n-1)}, \hat{Z}_s^{(n)} = Z_s^{(n)} - Z_s^{(n-1)}, n \geq 4$.

Then, by (7), we have

$$
\begin{aligned}
E\Big[\int_0^T \Big(\frac{\beta}{2}|\hat{Y}_s^{(n)}|^2 + |\hat{Z}_s^{(n)}|^2\Big) e^{\beta s}\, ds\Big]
&\le \frac{2}{\beta} E\Big[\int_0^T |f_2(s, Y_s^{(n)}, Z_s^{(n)}, Y_{s+\delta(s)}^{(n-1)}) \\
&\qquad\qquad - f_2(s, Y_s^{(n-1)}, Z_s^{(n-1)}, Y_{s+\delta(s)}^{(n-2)})|^2 e^{\beta s}\, ds\Big] \\
&\le \frac{6C^2}{\beta} E\Big[\int_0^T (|\hat{Y}_s^{(n)}|^2 + |\hat{Z}_s^{(n)}|^2) e^{\beta s}\, ds\Big] + \frac{6C^2 L}{\beta} E\Big[\int_0^T |\hat{Y}_s^{(n-1)}|^2 e^{\beta s}\, ds\Big].
\end{aligned}
$$

Set $\beta = 18C^2L + 18C^2 + 3$. Then

$$
\begin{aligned}
\frac{2}{3} E\Big[\int_0^T (|\hat{Y}_s^{(n)}|^2 + |\hat{Z}_s^{(n)}|^2) e^{\beta s}\, ds\Big]
&\le \frac{1}{3} E\Big[\int_0^T |\hat{Y}_s^{(n-1)}|^2 e^{\beta s}\, ds\Big] \le \frac{1}{3} E\Big[\int_0^T (|\hat{Y}_s^{(n-1)}|^2 + |\hat{Z}_s^{(n-1)}|^2) e^{\beta s}\, ds\Big].
\end{aligned}
$$

Hence,

$$
E\Big[\int_0^T (|\hat{Y}_s^{(n)}|^2 + |\hat{Z}_s^{(n)}|^2) e^{\beta s}\, ds\Big] \le \Big(\frac{1}{2}\Big)^{n-4} E\Big[\int_0^T (|\hat{Y}_s^{(4)}|^2 + |\hat{Z}_s^{(4)}|^2) e^{\beta s}\, ds\Big].
$$

It follows that $(Y_{\boldsymbol{\cdot}}^{(n)})_{n\ge 4}$ and $(Z_{\boldsymbol{\cdot}}^{(n)})_{n\ge 4}$ are respectively Cauchy sequences in $L^2_{\mathscr{F}}(0, T+K)$ and in $L^2_{\mathscr{F}}(0, T; \mathbb{R}^d)$. Denote their limits by $Y_{\boldsymbol{\cdot}}$ and $Z_{\boldsymbol{\cdot}}$, respectively. Since $L^2_{\mathscr{F}}(0, T+K)$ and $L^2_{\mathscr{F}}(0, T; \mathbb{R}^d)$ are both Banach spaces, we obtain $(Y_{\boldsymbol{\cdot}}, Z_{\boldsymbol{\cdot}}) \in L^2_{\mathscr{F}}(0, T+K) \times L^2_{\mathscr{F}}(0, T; \mathbb{R}^d)$. Note for all $t \in [0, T]$,

$$
\begin{aligned}
E\Big[\int_t^T |f_2(s, Y_s^{(n)}, Z_s^{(n)}, Y_{s+\delta(s)}^{(n-1)}) - f_2(s, Y_s, Z_s, Y_{s+\delta(s)})|^2 e^{\beta s}\, ds\Big] \\
\le 3C^2 E\Big[\int_t^T (|Y_s^{(n)} - Y_s|^2 + |Z_s^{(n)} - Z_s|^2 + L|Y_s^{(n-1)} - Y_s|^2) e^{\beta s}\, ds\Big] \to 0,
\end{aligned}
$$

when $n \to \infty$. Therefore, $(Y_{\boldsymbol{\cdot}}, Z_{\boldsymbol{\cdot}})$ satisfies the following anticipated BSDE:

$$
\begin{cases}
Y_t = \xi_T^{(2)} + \int_t^T f_2(s, Y_s, Z_s, Y_{s+\delta(s)})\, ds - \int_t^T Z_s\, dW_s, & 0 \le t \le T; \\
Y_t = \xi_t^{(2)}, & T \le t \le T+K.
\end{cases}
$$

By Theorem 4.2, we know

$$
Y_t = Y_t^{(2)}, \qquad \text{a.e., a.s.}
$$

Since $Y_t^{(1)} \ge Y_t^{(3)} \ge Y_t^{(4)} \ge Y_t$, it holds immediately

$$
Y_t^{(1)} \ge Y_t^{(2)}, \qquad \text{a.e., a.s.} \qquad \square
$$

If $f_2$ is nonincreasing in the anticipated term of $Y_\bullet$, Theorem 5.1 does not hold. The following example shows this.

EXAMPLE 5.2. Given $T > \delta > 0$, consider the following two anticipated BSDEs:

$$\begin{cases} Y_t = c + \int_t^T aE^{\mathscr{F}_s}[Y_{s+\delta}]\,ds - \int_t^T Z_s\,dW_s, & t \in [0,T]; \\ Y_t = c, & t \in [T, T+\delta], \end{cases} \tag{15}$$

and

$$\begin{cases} Y'_t = \int_t^T aE^{\mathscr{F}_s}[I_{[Y'_{s+\delta}<0]}Y'_{s+\delta}]\,ds - \int_t^T Z'_s\,dW_s, & t \in [0,T]; \\ Y'_t = 0, & t \in [T, T+\delta], \end{cases} \tag{16}$$

where $a = -\frac{2}{\delta}$, $c < 0$ are given constants. Obviously the solution to equation (16) is $(Y'_\bullet, Z'_\bullet) \equiv (0,0)$. When $t \in [T-\delta, T]$, equation (15) becomes

$$Y_t = c + \int_t^T ac\,ds - \int_t^T Z_s\,dW_s.$$

It is easy to see that $Y_t = c + ac(T-t)$, $Z_t \equiv 0$ is the solution of equation (15) when $t \in [T-\delta, T]$. But $Y_t > 0$ when $t \in [T-\delta, T-\delta/2)$.

If $f_2$ contains the anticipated term of $Z_\bullet$, Theorem 5.1 does not hold. This is shown in the following example.

EXAMPLE 5.3. Given $T > \delta > 0$, consider the two anticipated BSDEs

$$\begin{cases} Y_t = W_T^2 - T - \int_t^T \sqrt{\frac{\pi}{2\delta}} E^{\mathscr{F}_s}[|Z_{s+\delta} - Z_s|]\,ds - \int_t^T Z_s\,dW_s, \\ \qquad\qquad t \in [0,T]; \\ Y_t = W_t^2 - (T-t), \qquad t \in [T, T+\delta]; \\ Z_t = 2W_t, \qquad t \in [T, T+\delta], \end{cases} \tag{17}$$

and

$$\begin{cases} Y'_t = 4W_T^2 - \int_t^T \sqrt{\frac{\pi}{2\delta}} E^{\mathscr{F}_s}[|Z'_{s+\delta} - Z'_s|]\,ds - \int_t^T Z'_s\,dW_s, \\ \qquad\qquad t \in [0,T]; \\ Y'_t = 4W_t^2 - 4(T-t), \qquad t \in [T, T+\delta]; \\ Z'_t = 8W_t, \qquad t \in [T, T+\delta]. \end{cases} \tag{18}$$

We can check that the solution of (17) is $(Y_t, Z_t) = (W_t^2 - T - (T-t), 2W_t)$ and that the solution of (18) is $(Y'_t, Z'_t) = (4W_t^2 - 4(T-t), 8W_t)$. We have $Y_T < Y'_T$ but $Y_0 > Y'_0$.

THEOREM 5.4. *Under the assumptions of Theorem 5.1, if* $\xi_s^{(1)} \geq \xi_s^{(2)}$, $s \in [T, T+K]$ *and* $f_1(t, Y_t^{(1)}, Z_t^{(1)}, Y_{t+\delta(t)}^{(1)}) \geq f_2(t, Y_t^{(1)}, Z_t^{(1)}, Y_{t+\delta(t)}^{(1)})$, $t \in [0, T]$, *then*

$$Y_t^{(1)} \geq Y_t^{(2)}, \qquad a.e., \quad a.s.$$

*We also have a strict comparison. Given the assumptions of Theorem 5.1, suppose* $[T, T+K] \subset \{t + \delta(t), t \in [0, T]\}$ *and* $f_2$ *is strictly increasing in* $\theta$. *Then*

$$Y_0^{(1)} = Y_0^{(2)} \quad \Longleftrightarrow \quad \begin{cases} f_1(t, Y_t^{(1)}, Z_t^{(1)}, Y_{t+\delta(t)}^{(1)}) = f_2(t, Y_t^{(1)}, Z_t^{(1)}, Y_{t+\delta(t)}^{(1)}), \\ \qquad\qquad t \in [0, T], \\ \xi_s^{(1)} = \xi_s^{(2)}, \qquad s \in [T, T+K]. \end{cases}$$

PROOF. Set

$$\begin{cases} Y_t^{(3)} = \xi_T^{(2)} + \int_t^T f_2(s, Y_s^{(3)}, Z_s^{(3)}, Y_{s+\delta(s)}^{(1)})\, ds - \int_t^T Z_s^{(3)}\, dW_s, \\ \qquad\qquad t \in [0, T]; \\ Y_t^{(3)} = \xi_t^{(2)}, \qquad t \in [T, T+K]. \end{cases}$$

Set $\tilde{f}_t = f_1(t, Y_t^{(1)}, Z_t^{(1)}, Y_{t+\delta(t)}^{(1)}) - f_2(t, Y_t^{(1)}, Z_t^{(1)}, Y_{t+\delta(t)}^{(1)})$ and $y_\cdot = Y_\cdot^{(1)} - Y_\cdot^{(3)}$, $z_\cdot = Z_\cdot^{(1)} - Z_\cdot^{(3)}$, $\tilde{\xi}_\cdot = \xi_\cdot^{(1)} - \xi_\cdot^{(2)}$. Then the pair $(y_\cdot, z_\cdot)$ can be regarded as the solution to the linear BSDE

$$\begin{cases} y_t = \tilde{\xi}_T + \int_t^T (a_s y_s + b_s z_s + \tilde{f}_s)\, ds - \int_t^T z_s\, dW_s, \\ \qquad\qquad t \in [0, T]; \\ y_t = \tilde{\xi}_t, \qquad t \in [T, T+K], \end{cases}$$

where

$$a_s = \begin{cases} \dfrac{f_2(s, Y_s^{(1)}, Z_s^{(1)}, Y_{s+\delta(s)}^{(1)}) - f_2(s, Y_s^{(3)}, Z_s^{(1)}, Y_{s+\delta(s)}^{(1)})}{Y_s^{(1)} - Y_s^{(3)}}, & \text{if } Y_s^{(1)} \neq Y_s^{(3)}; \\ 0, & \text{if } Y_s^{(1)} = Y_s^{(3)}, \end{cases}$$

$$b_s = \begin{cases} \dfrac{f_2(s, Y_s^{(3)}, Z_s^{(1)}, Y_{s+\delta(s)}^{(1)}) - f_2(s, Y_s^{(3)}, Z_s^{(3)}, Y_{s+\delta(s)}^{(1)})}{Z_s^{(1)} - Z_s^{(3)}}, & \text{if } Z_s^{(1)} \neq Z_s^{(3)}; \\ 0, & \text{if } Z_s^{(1)} = Z_s^{(3)}. \end{cases}$$

Since $f_2$ satisfies (H1), $|a_s| \leq C$ and $|b_s| \leq C$. Set

$$X_t := \exp\left[\int_0^t b_s\, dW_s - \frac{1}{2}\int_0^t |b_s|^2\, ds + \int_0^t a_s\, ds\right] \geq 0.$$

We apply Itô's formula to $X_s y_s$ on $[t,T]$ and take conditional expectations on both sides:

$$y_t = E^{\mathscr{F}_t}\Big[\tilde{\xi}_T X_T + \int_t^T \tilde{f}_s X_s\,ds\Big].$$

Since $\tilde{\xi}_T \geq 0, \tilde{f}_t \geq 0$, a.e., a.s., we get $Y_t^{(1)} \geq Y_t^{(3)}$, a.e., a.s.

Then similarly to the proof of Theorem 5.1, we obtain

$$Y_t^{(1)} \geq Y_t^{(2)}, \qquad \text{a.e., a.s.}$$

Now we only need to prove the strict comparison theorem.

($\Longrightarrow$) Suppose $Y_0^{(1)} = Y_0^{(2)}$, by Lemma 3.3, we get

$$f_1(t, Y_t^{(1)}, Z_t^{(1)}, Y_{t+\delta(t)}^{(1)}) = f_2(t, Y_t^{(1)}, Z_t^{(1)}, Y_{t+\delta(t)}^{(2)}), \qquad t \in [0,T].$$

Since $Y_0^{(1)} \geq Y_0^{(3)} \geq Y_0^{(2)}$, we know $Y_0^{(1)} = Y_0^{(3)}$. Also by Lemma 3.3, we get

$$f_1(t, Y_t^{(1)}, Z_t^{(1)}, Y_{t+\delta(t)}^{(1)}) = f_2(t, Y_t^{(1)}, Z_t^{(1)}, Y_{t+\delta(t)}^{(1)}), \qquad t \in [0,T].$$

Therefore,

$$f_2(t, Y_t^{(1)}, Z_t^{(1)}, Y_{t+\delta(t)}^{(1)}) = f_2(t, Y_t^{(1)}, Z_t^{(1)}, Y_{t+\delta(t)}^{(2)}), \qquad t \in [0,T].$$

Note that for all $t \in [0,T], y \in \mathbb{R}, z \in \mathbb{R}^d$, $f_2(t,y,z,\cdot)$ is strictly increasing, hence, $Y_{t+\delta(t)}^{(1)} = Y_{t+\delta(t)}^{(2)}, t \in [0,T]$. In particular, $\xi_t^{(1)} = \xi_t^{(2)}, t \in [T, T+K]$.

($\Longleftarrow$) Suppose $f_1(t, Y_t^{(1)}, Z_t^{(1)}, Y_{t+\delta(t)}^{(1)}) = f_2(t, Y_t^{(1)}, Z_t^{(1)}, Y_{t+\delta(t)}^{(1)}), t \in [0,T]$ and $\xi_s^{(1)} = \xi_s^{(2)}, s \in [T, T+K]$. Then

$$y_t = Y_t^{(1)} - Y_t^{(3)} = E^{\mathscr{F}_t}\Big[\tilde{\xi}_T X_T + \int_t^T \tilde{f}_s X_s\,ds\Big] \equiv 0.$$

Therefore,

$$\begin{cases} Y_t^{(1)} = \xi_T^{(2)} + \displaystyle\int_t^T f_2(s, Y_s^{(1)}, Z_s^{(3)}, Y_{s+\delta(s)}^{(1)})\,ds - \int_t^T Z_s^{(3)}\,dW_s, \\ \qquad\qquad t \in [0,T]; \\ Y_t^{(1)} = \xi_t^{(2)}, \qquad t \in [T, T+K]. \end{cases}$$

By Theorem 4.2, $Y_t^{(1)} = Y_t^{(2)}$, a.e., a.s., in particular, $Y_0^{(1)} = Y_0^{(2)}$. □

Corollary 5.5. *Let $(Y_\cdot^{(1)}, Z_\cdot^{(1)})$ and $(Y_\cdot^{(2)}, Z_\cdot^{(2)})$ be respectively the solutions for the following two 1-dimensional anticipated BSDEs:*

$$\begin{cases} -dY_t^{(j)} = f(t, Y_t^{(j)}, Z_t^{(j)}, Y_{t+\delta_j(t)}^{(j)})\,dt - Z_t^{(j)}\,dW_t, & 0 \leq t \leq T; \\ Y_t^{(j)} = \xi_t, & T \leq t \leq T+K, \end{cases}$$

*where* $j=1,2$. *Suppose* $\xi_\cdot$ *is in* $S^2_{\mathscr{F}}(T,T+K)$, $f$ *satisfies* (H1) *and* (H2), *for all* $t\in[0,T], y\in\mathbb{R}, z\in\mathbb{R}^d$, $f(t,y,z,\cdot)$ *is increasing, and* $\delta_1,\delta_2$ *satisfy* (i) *and* (ii). *If* $Y^{(1)}_{t+\delta_1(t)}\geq Y^{(1)}_{t+\delta_2(t)}$, *a.e., a.s., then*

$$Y^{(1)}_t\geq Y^{(2)}_t,\qquad a.e.,\quad a.s.$$

PROOF. Set

$$\begin{cases} Y^{(3)}_t=\xi_T+\displaystyle\int_t^T f(s,Y^{(3)}_s,Z^{(3)}_s,Y^{(1)}_{s+\delta_2(s)})\,ds-\int_t^T Z^{(3)}_s\,dW_s,\\ \qquad\qquad t\in[0,T];\\ Y^{(3)}_t=\xi_t,\qquad t\in[T,T+K].\end{cases}$$

From Lemma 3.2, there exists a unique pair of $\mathscr{F}_t$-adapted processes $(Y^{(3)}_\cdot, Z^{(3)}_\cdot)\in S^2_{\mathscr{F}}(0,T)\times L^2_{\mathscr{F}}(0,T;\mathbb{R}^{1\times d})$ that satisfies the above BSDE. Since $f(s,y,z,Y^{(1)}_{s+\delta_1(s)})\geq f(s,y,z,Y^{(1)}_{s+\delta_2(s)})$, by Lemma 3.4, we know

$$Y^{(1)}_t\geq Y^{(3)}_t,\qquad \text{a.e.,}\quad \text{a.s.}$$

The remaining proof is similar to Theorem 5.1, we omit it. □

**6. Stochastic control problems.** El Karoui, Peng and Quenez [7] applied the duality between SDEs and BSDEs to stochastic control problems. Now we consider if it is feasible to use the duality between SDDEs and anticipated BSDEs to solve these problems. Let $\theta>0$ be a given constant. Now we consider the following stochastic control problem: the laws of the controlled process belong to a family of equivalent measures whose densities are

$$\begin{cases} dX^u_s=(\alpha(s,u_s)X^u_s+b(s-\theta,u_{s-\theta})X^u_{s-\theta})\,ds+X^u_s\sigma^T(s,u_s)\,dW_s,\\ \qquad\qquad s\in[t,T+\theta];\\ X^u_t=1,\\ X^u_s=0,\qquad s\in[t-\theta,t),\end{cases}$$

where the coefficients $\alpha(s,u):\mathbb{R}\times\mathbb{R}^k\longrightarrow\mathbb{R}, b(s,u):\mathbb{R}\times\mathbb{R}^k\longrightarrow\mathbb{R}^+$ and $\sigma(s,u):\mathbb{R}\times\mathbb{R}^k\longrightarrow\mathbb{R}^{d\times 1}$ are adapted processes uniformly continuous with respect to $(s,u)$. A feasible control $(u_s, s\in[-\theta,T+\theta])$ is a continuous adapted process valued in a compact subset $U$ in $\mathbb{R}^k$. The set of feasible controls is denoted by $\mathcal{U}$. The problem is to maximize over all feasible control processes $u$ the objective function

$$J(u)=E\biggl[X^u_TQ(T)+\int_T^{T+\theta}X^u_{s-\theta}Q(s)b(s-\theta,u_{s-\theta})\,ds+\int_0^T X^u_s l(s,u_s)\,ds\biggr],$$

where $Q(\cdot) \in S^2_{\mathscr{F}}(T, T+\theta)$ is the terminal condition, $(l(\omega, s, u_s), s \in [0,T])$ is the running cost associated with the control process $u$ and $l(s,u)$ is an adapted process uniformly continuous with respect to $(s,u)$. Assume $\alpha(s,u), b(s,u), |\sigma(s,u)|$ and $l(s,u)$ are uniformly bounded by $\mu$. Notice that, by Theorem 2.1, $J(u) = Y_0^u$, where $(Y_\cdot^u, Z_\cdot^u)$ is the solution to the following linear anticipated BSDE:

$$\begin{cases} -dY_t^u = f^u(t, Y_t^u, Z_t^u, Y_{t+\theta}^u)\,dt - Z_t^u\,dW_t, & t \in [0,T]; \\ Y_t^u = Q(t), & t \in [T, T+\theta], \end{cases}$$

where $f^u(t,y,z,\eta_r) = \alpha(t,u_t)y + z\sigma(t,u_t) + b(t,u_t)E^{\mathscr{F}_t}[\eta_r] + l(t,u_t), t \in [0,T], y \in \mathbb{R}, z \in \mathbb{R}^d, \eta_\cdot \in L^2_{\mathscr{F}}(t, T+\theta), r \in [t, T+\theta]$ and

$$Y_t^u = E^{\mathscr{F}_t}\biggl[X_T^u Q(T) + \int_T^{T+\theta} X_{s-\theta}^u Q(s) b(s-\theta, u_{s-\theta})\,ds + \int_t^T X_s^u l(s, u_s)\,ds\biggr].$$

THEOREM 6.1. *Set* $f(t,y,z,\eta_r) = \operatorname{esssup}\{f^u(t,y,z,\eta_r), u \in \mathcal{U}\}, t \in [0,T]$, $y \in \mathbb{R}, z \in \mathbb{R}^d, \eta_\cdot \in L^2_{\mathscr{F}}(t, T+\theta), r \in [t, T+\theta]$. *Then anticipated BSDE*

$$\begin{cases} -dY_t = f(t, Y_t, Z_t, Y_{t+\theta})\,dt - Z_t\,dW_t, & t \in [0,T]; \\ Y_t = Q(t), & t \in [T, T+\theta], \end{cases} \tag{19}$$

*has a unique solution* $(Y_\cdot, Z_\cdot)$. *Moreover,* $Y_\cdot$ *is the value function* $Y_\cdot^*$ *of the control problem, that is, for each* $t \in [0,T]$,

$$Y_t = Y_t^* = \operatorname{esssup}\{Y_t^u, u \in \mathcal{U}\}.$$

PROOF. On one hand, since $\alpha, b, |\sigma|$ and $l$ are uniformly bounded by $\mu$, for all $t \in [0,T], s \in [T, T+\theta], y, y' \in \mathbb{R}, z, z' \in \mathbb{R}^d, \eta_\cdot, \eta'_\cdot \in L^2_{\mathscr{F}}(t, T+\theta)$, and $r \in [t, T+\theta]$,

$$\begin{aligned} f(t,y,z,\eta_r) - f(t,y',z',\eta'_r) &\le \operatorname{esssup}\{\alpha(t,u_t)(y-y') + (z-z')\sigma(t,u_t) \\ &\qquad + b(t,u_t)E^{\mathscr{F}_t}[\eta_r - \eta'_r], u \in \mathcal{U}\} \\ &\le \mu(|y-y'| + |z-z'| + E^{\mathscr{F}_t}[|\eta_r - \eta'_r|]). \end{aligned}$$

Notice $E[\int_0^T |f(t,0,0,0)|^2\,dt] \le \mu^2 T$, then by Theorem 4.2, the anticipated BSDE (19) has a unique solution $(Y., Z.) \in S^2_{\mathscr{F}}(0, T+K) \times L^2_{\mathscr{F}}(0, T; \mathbb{R}^d)$.

Since for all $u \in \mathcal{U}, f^u(t,y,z,\eta) \le f(t,y,z,\eta)$ and $f^u(t,y,z,\eta)$ is increasing in $\eta$, by Theorem 5.1, we get $Y_t \ge Y_t^u$, a.e., a.s. Thus, $Y_t \ge Y_t^*$, a.e., a.s.

On the other hand, by the definition of $f$, we know for all $\varepsilon > 0$, for each $(\omega, t) \in \Omega \times [0, T)$,

$$\begin{aligned}\{u \in \mathcal{U}; f(\omega, t, Y_t(\omega), Z_t(\omega), Y_{t+\theta}(\omega)) \\ \leq \alpha(t,u)Y_t(\omega) + Z_t(\omega)\sigma(t,u) \\ + b(t,u)E^{\mathscr{F}_t}[Y_{t+\theta}(\omega)] + l(\omega, t, u) + \varepsilon\} \neq \varnothing.\end{aligned}$$

Then by a Measurable Selection Theorem, for example, that can be found in Dellacherie [5] or in Beneš [1, 2], there exists a $u^\varepsilon \in \mathcal{U}$ such that

$$f(t, Y_t, Z_t, Y_{t+\theta}) \leq f^{u^\varepsilon}(t, Y_t, Z_t, Y_{t+\theta}) + \varepsilon, \qquad \text{a.e., a.s.}$$

Denote the solution to the anticipated BSDE corresponding to $(f^{u^\varepsilon}, Q_\bullet)$ by $(Y_\bullet^{u^\varepsilon}, Z_\bullet^{u^\varepsilon})$.

First, consider the case when $t \in [T-\theta, T]$. Thus, $t+\theta \in [T, T+\theta]$, $Y_{t+\theta} = Y_{t+\theta}^{u^\varepsilon}$ and

$$\begin{aligned}f^{u^\varepsilon}&(t, Y_t^{u^\varepsilon}, Z_t^{u^\varepsilon}, Y_{t+\theta}^{u^\varepsilon}) - f(t, Y_t, Z_t, Y_{t+\theta}) \\ &\geq f^{u^\varepsilon}(t, Y_t^{u^\varepsilon}, Z_t^{u^\varepsilon}, Y_{t+\theta}^{u^\varepsilon}) - f^{u^\varepsilon}(t, Y_t, Z_t, Y_{t+\theta}) - \varepsilon \\ &= f^{u^\varepsilon}(t, Y_t^{u^\varepsilon}, Z_t^{u^\varepsilon}, Y_{t+\theta}^{u^\varepsilon}) - f^{u^\varepsilon}(t, Y_t, Z_t, Y_{t+\theta}^{u^\varepsilon}) - \varepsilon \\ &= g_t^{(1)}(Y_t^{u^\varepsilon} - Y_t) + g_t^{(2)}(Z_t^{u^\varepsilon} - Z_t) - \varepsilon,\end{aligned}$$

where, for $t \in [T-\theta, T]$,

$$g_t^{(1)} = \begin{cases} \dfrac{f^{u^\varepsilon}(t, Y_t^{u^\varepsilon}, Z_t^{u^\varepsilon}, Y_{t+\theta}^{u^\varepsilon}) - f^{u^\varepsilon}(t, Y_t, Z_t^{u^\varepsilon}, Y_{t+\theta}^{u^\varepsilon})}{Y_t^{u^\varepsilon} - Y_t}, & \text{if } Y_t^{u^\varepsilon} \neq Y_t; \\ 0, & \text{if } Y_t^{u^\varepsilon} = Y_t, \end{cases}$$

$$g_t^{(2)} = \begin{cases} \dfrac{f^{u^\varepsilon}(t, Y_t, Z_t^{u^\varepsilon}, Y_{t+\theta}^{u^\varepsilon}) - f^{u^\varepsilon}(t, Y_t, Z_t, Y_{t+\theta}^{u^\varepsilon})}{Z_t^{u^\varepsilon} - Z_t}, & \text{if } Z_t^{u^\varepsilon} \neq Z_t; \\ 0, & \text{if } Z_t^{u^\varepsilon} = Z_t. \end{cases}$$

That is, for $t \in [T-\theta, T]$,

$$\begin{aligned}Y_t^{u^\varepsilon} - Y_t \geq &\int_t^T (g_s^{(1)}(Y_s^{u^\varepsilon} - Y_s) + g_s^{(2)}(Z_s^{u^\varepsilon} - Z_s) - \varepsilon)\, ds \\ &- \int_t^T (Z_s^{u^\varepsilon} - Z_s)\, dW_s.\end{aligned}$$

Hence, $Y_t^{u^\varepsilon} - Y_t \geq \tilde{Y}_t^{(1)}$, where $\tilde{Y}_t^{(1)}$ is the solution of BSDE:

$$\tilde{Y}_t^{(1)} = \int_t^T (g_s^{(1)}\tilde{Y}_s^{(1)} + g_s^{(2)}\tilde{Z}_s^{(1)} - \varepsilon)\, ds - \int_t^T \tilde{Z}_s^{(1)}\, dW_s, \qquad t \in [T-\theta, T].$$

Since $|g_t^{(1)}| \le \mu$, $|g_t^{(2)}| \le \mu$, we get

$$\tilde{Y}_t^{(1)} = -\varepsilon E^{\mathscr{F}_t}\biggl[\int_t^T \tilde{X}_s^{(1)}\,ds\biggr], \qquad t \in [T-\theta, T],$$

where

$$\tilde{X}_t^{(1)} = \exp\biggl[\int_0^t g_s^{(2)}\,dW_s - \frac{1}{2}\int_0^t |g_s^{(2)}|^2\,ds + \int_0^t g_s^{(1)}\,ds\biggr] \ge 0.$$

Therefore, there exists a constant $\rho_1 > 0$ depending only on $\mu$, $\theta$ and $T$ such that

$$Y_t^{u^\varepsilon} - Y_t \ge \tilde{Y}_t^{(1)} \ge -\rho_1\varepsilon, \qquad t \in [T-\theta, T].$$

Second, consider the case when $t \in [T-2\theta, T-\theta]$. Then $t+\theta \in [T-\theta, T]$, $Y_{t+\theta} \le Y_{t+\theta}^{u^\varepsilon} + \rho_1\varepsilon$. Since for all $t \in [0,T], y \in \mathbb{R}, z \in \mathbb{R}^d$, $f^u(t,y,z,\cdot)$ is a increasing and linear function, we have

$$\begin{aligned}
&f^{u^\varepsilon}(t, Y_t^{u^\varepsilon}, Z_t^{u^\varepsilon}, Y_{t+\theta}^{u^\varepsilon}) - f(t, Y_t, Z_t, Y_{t+\theta})\\
&\quad \ge f^{u^\varepsilon}(t, Y_t^{u^\varepsilon}, Z_t^{u^\varepsilon}, Y_{t+\theta}^{u^\varepsilon}) - f^{u^\varepsilon}(t, Y_t, Z_t, Y_{t+\theta}) - \varepsilon\\
&\quad \ge f^{u^\varepsilon}(t, Y_t^{u^\varepsilon}, Z_t^{u^\varepsilon}, Y_{t+\theta}^{u^\varepsilon}) - f^{u^\varepsilon}(t, Y_t, Z_t, Y_{t+\theta}^{u^\varepsilon} + \rho_1\varepsilon) - \varepsilon\\
&\quad \ge f^{u^\varepsilon}(t, Y_t^{u^\varepsilon}, Z_t^{u^\varepsilon}, Y_{t+\theta}^{u^\varepsilon}) - f^{u^\varepsilon}(t, Y_t, Z_t, Y_{t+\theta}^{u^\varepsilon}) - \mu\rho_1\varepsilon - \varepsilon\\
&\quad = g_t^{(1)}(Y_t^{u^\varepsilon} - Y_t) + g_t^{(2)}(Z_t^{u^\varepsilon} - Z_t) - (\mu\rho_1 + 1)\varepsilon,
\end{aligned}$$

where, for $t \in [T-2\theta, T-\theta]$,

$$g_t^{(1)} = \begin{cases} \dfrac{f^{u^\varepsilon}(t, Y_t^{u^\varepsilon}, Z_t^{u^\varepsilon}, Y_{t+\theta}^{u^\varepsilon}) - f^{u^\varepsilon}(t, Y_t, Z_t^{u^\varepsilon}, Y_{t+\theta}^{u^\varepsilon})}{Y_t^{u^\varepsilon} - Y_t}, & \text{if } Y_t^{u^\varepsilon} \ne Y_t;\\ 0, & \text{if } Y_t^{u^\varepsilon} = Y_t, \end{cases}$$

$$g_t^{(2)} = \begin{cases} \dfrac{f^{u^\varepsilon}(t, Y_t, Z_t^{u^\varepsilon}, Y_{t+\theta}^{u^\varepsilon}) - f^{u^\varepsilon}(t, Y_t, Z_t, Y_{t+\theta}^{u^\varepsilon})}{Z_t^{u^\varepsilon} - Z_t}, & \text{if } Z_t^{u^\varepsilon} \ne Z_t;\\ 0, & \text{if } Z_t^{u^\varepsilon} = Z_t. \end{cases}$$

Therefore, for $t \in [T-2\theta, T-\theta]$,

$$\begin{aligned}
Y_t^{u^\varepsilon} - Y_t &\ge Y_{T-\theta}^{u^\varepsilon} - Y_{T-\theta} - \int_t^{T-\theta} (Z_s^{u^\varepsilon} - Z_s)\,dW_s\\
&\qquad + \int_t^{T-\theta} (g_s^{(1)}(Y_s^{u^\varepsilon} - Y_s) + g_s^{(2)}(Z_s^{u^\varepsilon} - Z_s) - (\mu\rho_1 + 1)\varepsilon)\,ds.
\end{aligned}$$

Hence, $Y_t^{u^\varepsilon} - Y_t \ge \tilde{Y}_t^{(2)}$, where $\tilde{Y}_t^{(2)}$ is the solution to the following BSDE: For $t \in [T-2\theta, T-\theta]$,

$$\begin{aligned}
\tilde{Y}_t^{(2)} &= Y_{T-\theta}^{u^\varepsilon} - Y_{T-\theta}\\
&\qquad + \int_t^{T-\theta} (g_s^{(1)}\tilde{Y}_s^{(2)} + g_s^{(2)}\tilde{Z}_s^{(2)} - (\mu\rho_1 + 1)\varepsilon)\,ds - \int_t^{T-\theta} \tilde{Z}_s^{(2)}\,dW_s.
\end{aligned}$$

Note $|g_t^{(1)}| \le \mu$ and $|g_t^{(2)}| \le \mu$. We have, for all $t \in [T-2\theta, T-\theta]$,

$$\tilde{Y}_t^{(2)} = E^{\mathscr{F}_t}\left[(Y_{T-\theta}^{u^\varepsilon} - Y_{T-\theta})\tilde{X}_{T-\theta}^{(2)} - \int_t^{T-\theta} (\mu\rho_1 + 1)\varepsilon \tilde{X}_s^{(2)}\, ds\right],$$

where

$$\tilde{X}_t^{(2)} = \exp\left[\int_0^t g_s^{(2)}\, dW_s - \frac{1}{2}\int_0^t |g_s^{(2)}|^2\, ds + \int_0^t g_s^{(1)}\, ds\right] \ge 0.$$

Since $Y_{T-\theta}^{u^\varepsilon} - Y_{T-\theta} \ge -\rho_1\varepsilon$, there exists a constant $\rho_2 > 0$ depending only on $\mu$, $\theta$ and $T$ such that

$$Y_t^{u^\varepsilon} - Y_t \ge \tilde{Y}_t^{(2)} \ge -\rho_2\varepsilon, \qquad t \in [T-2\theta, T-\theta].$$

Similarly, we get constants $\rho_3, \rho_4, \dots, \rho_{[\frac{T}{\theta}]+1} > 0$ such that

$$Y_t^{u^\varepsilon} - Y_t \ge -\rho_n\varepsilon, \qquad t \in [T-n\theta, T-(n-1)\theta], n = 3, 4, \dots, \left[\frac{T}{\theta}\right];$$

$$Y_t^{u^\varepsilon} - Y_t \ge -\rho_{[T/\theta]+1}\varepsilon, \qquad t \in \left[0, T - \left[\frac{T}{\theta}\right]\theta\right].$$

Setting $\rho = \max\{\rho_2, \rho_3, \dots, \rho_{[T/\theta]+1}\}$, we obtain

$$Y_t^{u^\varepsilon} - Y_t \ge -\rho\varepsilon, \qquad t \in [0, T].$$

Since $Y_t^{u^\varepsilon} \le Y_t$, a.e., a.s., setting $\varepsilon \to 0$, we get $Y_t^{u^\varepsilon} \to Y_t$, a.e., a.s. Thus,

$$Y_t = Y_t^*, \qquad \text{a.e., a.s.} \qquad \square$$

**Acknowledgments.** The authors gratefully acknowledge the help from Doctor Kessel Cathy to enhance the readability of this paper. The second author wishes to thank Professor Xuerong Mao, Doctor Chenggui Yuan and Professor Huaizhong Zhao for their enthusiastic hospitality and discussions during her visits.

## REFERENCES

[1] Beneš, V. E. (1970). Existence of optimal strategies based on specified information, for a class of stochastic decision problems. *SIAM J. Control Optim.* **8** 179–188. MR0265043

[2] Beneš, V. E. (1971). Existence of optimal stochastic control laws. *SIAM J. Control Optim.* **9** 446–472. MR0300726

[3] Cao, Z. G. and Yan, J.-A. (1999). A comparison theorem for solutions of backward stochastic differential equations. *Adv. Math. (China)* **28** 304–308. MR1767637

[4] Cvitanić, J. and Karatzas, I. (1996). Backward stochastic differential equations with reflection and Dynkin games. *Ann. Probab.* **24** 2024–2056. MR1415239

[5] Dellacherie, C. (1972). *Capacités et Processus Stochastiques. Ergebnisse der Mathematik und ihrer Grenzgebiete* **67**. Springer, Berlin. MR0448504

[6] El Karoui, N., Kapoudjian, C., Pardoux, E., Peng, S. and Quenez, M. C. (1997). Reflected solutions of backward SDE's, and related obstacle problems for PDE's. *Ann. Probab.* **25** 702–737. MR1434123

[7] El Karoui, N., Peng, S. and Quenez, M. C. (1997). Backward stochastic differential equations in finance. *Math. Finance* **7** 1–71. MR1434407

[8] Hu, Y. and Peng, S. (2006). On the comparison theorem for multidimensional BSDEs. *C. R. Math. Acad. Sci. Paris Ser. I* **343** 135–140. MR2243308

[9] Lepeltier, J.-P. and Martín, J. S. (2004). Backward SDEs with two barriers and continuous coefficient: An existence result. *J. Appl. Probab.* **41** 162–175. MR2036279

[10] Lin, Q. Q. (2001). A comparison theorem for backward stochastic differential equations. *J. Huazhong Univ. Sci. Tech.* **29** 1–3. MR1858002

[11] Liu, J. and Ren, J. (2002). Comparison theorem for solutions of backward stochastic differential equations with continuous coefficient. *Statist. Probab. Lett.* **56** 93–100. MR1881535

[12] Peng, S. (1999). Monotonic limit theorem of BSDE and nonlinear decomposition theorem of Doob–Meyer's type. *Probab. Theory Related Fields* **113** 473–499. MR1717527

[13] Peng, S. (2004). Nonlinear expectations, nonlinear evaluations and risk measures. In *Stochastic Methods in Finance. Lecture Notes in Math.* **1856** 165–253. Springer, Berlin. MR2113723

[14] Peng, S. and Xu, M. Y. (2005). The smallest $g$-supermartingale and reflected BSDE with single and double $L^2$ obstacles. *Ann. Inst. H. Poincaré Probab. Statist.* **41** 605–630. MR2139035

[15] Situ, R. (1999). Comparison theorem of solutions to BSDE with jumps, and viscosity solution to a generalized Hamilton–Jacobi–Bellman equation. In *Control of Distributed Parameter and Stochastic Systems (Hangzhou, 1998)* 275–282. Kluwer Academic, Boston, MA. MR1777420

[16] Zhang, T. S. (2003). A comparison theorem for solutions of backward stochastic differential equations with two reflecting barriers and its applications. In *Probabilistic Methods in Fluids* 324–331. World Science, River Edge, NJ. MR2083381

School of Mathematics
Shandong University
Jinan, 250100
China
E-mail: peng@sdu.edu.cn

School of Mathematics
Shandong University
Jinan, 250100
China
And
Present Address:
Department of Engineering
Cambridge University
Trumpington Street
Cambridge CB2 1PZ
United Kingdom
E-mail: yangzhezhe@gmail.com